\RequirePackage{fix-cm}
\documentclass[smallextended]{svjour3}       
\smartqed  
\usepackage{graphicx}
\usepackage{amsmath}
\usepackage{hyperref}
\newcommand{\ampl}{{\sc ampl}}
\newcommand{\loqo}{{\sc loqo}}

\newcommand{\ds}{\displaystyle}
\newcommand{\dxi}{{\Delta \xi}}
\newcommand{\deta}{{\Delta \eta}}
\newcommand{\dx}{{\Delta x}}
\newcommand{\dy}{{\Delta y}}

\newenvironment{alignedt}{
                \setlength{\arraycolsep}{0.1em}
                \begin{array}[t]{rcllll}
                \displaystyle
        }{
                \end{array}
        }
\begin{document}

\title{Fast Fourier Optimization\thanks{The author was supported by a grant from
	NASA.}
}
\subtitle{Sparsity Matters}


\author{Robert J. Vanderbei
}


\institute{Robert J. Vanderbei \at
	Department of Ops. Res. and Fin. Eng.,
	Princeton University,
	Princeton, NJ 08544. \\
              Tel.: +609-258-2345 \\
              \email{rvdb@princeton.edu}           
}

\date{Received: date / Accepted: date}

\maketitle

\begin{abstract}
Many interesting and fundamentally practical optimization problems, 
     ranging from optics, to signal
processing, to radar and acoustics, involve constraints on the Fourier transform of 
a function.  It is well-known that the {\em fast Fourier transform} (fft) is a recursive
algorithm that can dramatically improve the efficiency for computing the discrete
Fourier transform.  However, because it is recursive, it is difficult
to embed into a linear optimization problem.  In this paper, we explain the
main idea behind the fast Fourier transform and show how to adapt it
in such a manner as to make it encodable as constraints in an
optimization problem.   We demonstrate a real-world problem from the field of 
high-contrast imaging.  On this problem, dramatic improvements are translated to an ability to
solve problems with a much finer grid of discretized points.
As we shall show, in general, the ``fast Fourier'' version of the optimization constraints
produces a larger but sparser constraint matrix and therefore one
can think of the fast Fourier transform as a method of sparsifying the
constraints in an optimization problem, which is usually a good thing.

\end{abstract}

\section{Fourier Transforms in Engineering}

Many problems in engineering involve maximizing (or minimizing) a linear
functional of an unknown real-valued design function $f$ subject to constraints on 
its Fourier transform $\widehat{f}$ at certain points in transform space
(\cite{BN01}).  
Examples include antenna array synthesis (see, e.g.,
\cite{LB97,ref:Mai05,SC07}), FIR filter design (see, e.g.,
	\cite{CS99,ref:WBV96,ref:WBV99}), and
coronagraph design (see, e.g.,
\cite{ref:Indebetouw,SK01,KVSL02,KVLS04,Sou05,TEANK06,GPKCR06,TT07,MDC09}).
If the design function 
$f$ can be constrained to vanish outside a compact interval $C = (-a,a)$ of the real line
centered at the origin, then we can write the Fourier transform as
\[
    \widehat{f}(\xi) = \int_{-a}^a e^{2 \pi i x \xi} f(x) dx
\]
and an optimization problem might look like
\begin{equation} \label{2}
	\begin{array}{ll}
	    \text{maximize }   & \int_{-a}^a c(x) f(x) dx \\
	    \text{subject to } & \begin{alignedt}
		-\varepsilon & \le \Re \widehat{f}(\xi) & \le \varepsilon,  \qquad \xi \in D \\
		-\varepsilon & \le \Im \widehat{f}(\xi) & \le \varepsilon,  \qquad \xi \in D \\
				  0 & \le f(x) & \le 1, \qquad x \in C ,
				 \end{alignedt} 
	\end{array}
\end{equation}
where $D$ is a given subset of the real line, $\varepsilon$ is a given constant,
and $\Re(z)$ and $\Im(z)$ denote the real and imaginary parts of the complex
number $z$.
In Section \ref{sec_appl}, we will discuss a specific
real-world problem that fits a two-dimensional version of
this optimization paradigm and for which dramatic
computational improvements can be made.

Problem \eqref{2} is linear but it is infinite dimensional.   The first step to
making a tractable problem is to discretize both sets $C$ and $D$ so that the
continuous Fourier transform can be approximated by a discrete Riemann sum:
\begin{equation} \label{1}
    \widehat{f}_j = \sum_{k=-n}^n e^{2 \pi i k \dx j \dxi} f_k \dx, 
    		\qquad -n \le j \le n.
\end{equation}
Here, $n$ denotes the level of discretization,
\[
    \dx  = \frac{2a}{2n+1}, 
\]
$\dxi$ denotes the discretization spacing in transform space,
$f_k = f(k \dx)$, and $\widehat{f}_j \approx \widehat{f}(j \dxi)$.

Computing the discrete approximation \eqref{1} by simply summing the terms in its definition
requires on the order of $N^2$ operations, where $N = 2n+1$ is the number of
discrete points in both the function space and the transform space (later we
will generalize to allow a different number of points in the
discretization of $C$ and $D$).

Choosing $\dxi$ too large creates redundancy in the discrete approximation due
to periodicity of the complex exponential function and hence one generally
chooses $\dxi$ such that
\[
    \dx \dxi \le \frac{1}{N}.
\]
In many real-world applications, $\dxi$ is chosen so that this inequality is an
equality:  $\dxi = 1/(N \dx)$.
In this case, the Riemann sum approximation is called the {\em discrete Fourier
transform}.

\section{A Fast Fourier Transform}

Over the past half century there has been an explosion of research into
algorithms for efficiently computing Fourier transforms.  Any algorithm that
can do the job in a constant times $N \log N$ multiplications/additions is
called a {\em fast Fourier transform} (see, e.g., \cite{CT65,RM67,Pap79,DV90}).   
There are several algorithms that can
be called fast Fourier transforms.   Here, we present one that applies naturally
to Fourier transforms expressed as in \eqref{1}.  In this section, we assume
that $\dxi = 1/(N \dx)$.

A sum from $-n$ to $n$ has an odd number of terms: $N = 2n+1$.   Suppose, for
this section, that $N$ is a power of three:
\[
    N = 3^m.
\]
Fast Fourier transform algorithms assume that it is possible to factor 
$N$ into a product 
\[
    N = N_0 N_1.
\]
For the algorithm of this section, we put
\[
    N_0 = 3, \qquad \text{and} \qquad N_1 = 3^{m-1} .
\]
The first key idea in fast Fourier transform algorithms is to write the single sum
\eqref{2}
as a double sum and simultaneously to represent the discrete set of transform
values as a two-dimensional array of values rather than as a one-dimensional
vector.   Specifically, we decompose $k$ as 
\[
    k = N_0 k_1 + k_0
\]
so that 
\[
    -n \le k \le n 
    \quad \Longleftrightarrow \quad 
    -n_0 \le k_0 \le n_0 \quad \text{ and } \quad -n_1 \le k_1 \le n_1,
\]
where 
\[
    n_0 = (N_0-1)/2 = (3-1)/2 = 1
\]
and
\[
    n_1 = (N_1-1)/2 = (3^{m-1}-1)/2 .
\]
Similarly, we decompose $j$ as 
\[
    j = N_1 j_1 + j_0 
\]
so that
\[
    -n \le j \le n 
    \quad \Longleftrightarrow \quad 
    -n_1 \le j_0 \le n_1 \quad \text{ and } \quad -1 \le j_1 \le 1.
\]
With these notations, we rewrite the Fourier transform \eqref{1} as a double sum:
\begin{equation} \label{3}
    \widehat{f}_{j_0,j_1} 
	=
	\sum_{k_0=-1}^1 
	\sum_{k_1=-n_1}^{n_1} 
	e^{2 \pi i (N_0 k_1 + k_0) \dx (N_1 j_1 + j_0) \dxi} 
	f_{k_0,k_1} \dx ,
\end{equation}
where $f_{k_0,k_1} = f_{N_0 k_1 + k_0}$ and 
$\widehat{f}_{j_0,j_1} = \widehat{f}_{ N_1j_1 + j_0 }$.
Distributing the multiplications over the sums, we can rewrite the exponential
as
\vspace*{0.2in}
\begin{eqnarray*}
	&& e^{2 \pi i (N_0 k_1 + k_0) \dx (N_1 j_1 + j_0) \dxi} \\
	&& \qquad \qquad = 
	e^{2 \pi i N_0 k_1 \dx (N_1 j_1 + j_0) \dxi}  \;
	e^{2 \pi i k_0 \dx (N_1 j_1 + j_0) \dxi} 
	\\
	&& \qquad \qquad = 
	e^{2 \pi i N_0 k_1 \dx N_1 j_1 \dxi}  \;
	e^{2 \pi i N_0 k_1 \dx j_0 \dxi}  \;
	e^{2 \pi i k_0 \dx (N_1 j_1 + j_0) \dxi} 
	\\
	&& \qquad \qquad = 
	e^{2 \pi i N_0 k_1 \dx j_0 \dxi}  \;
	e^{2 \pi i k_0 \dx (N_1 j_1 + j_0) \dxi} ,
	\\ ~
\end{eqnarray*}
where the last equality follows from our assumption that 
$N_0 N_1 \dx \dxi = N \dx \dxi = 1$.
Substituting into \eqref{3}, we get
\[
    \widehat{f}_{j_0,j_1} 
	=
	\sum_{k_0=-1}^1 
	e^{2 \pi i k_0 \dx (N_1 j_1 + j_0) \dxi} \;
	\left(
	\sum_{k_1=-n_1}^{n_1} 
	e^{2 \pi i N_0 k_1 \dx j_0 \dxi} 
	\cdot 
	f_{k_0,k_1} 
	\right) \dx .
\]
We can compute this nested sum in two steps:
\begin{equation} \label{5}
  \begin{alignedt}
    g_{j_0, k_0} 
	& \; = \; &
	\ds \sum_{k_1=-n_1}^{n_1} 
	e^{2 \pi i N_0 k_1 \dx j_0 \dxi}  \;
	f_{k_0,k_1} \dx , 
	& \quad \begin{array}{c} -n_1 \le j_0 \le n_1, \\ -1 \le k_0 \le 1 \end{array}
	\\[0.2in]
    \widehat{f}_{j_0,j_1} 
	& \; = \; &
	\ds \sum_{k_0=-1}^1 
	e^{2 \pi i k_0 \dx j \dxi} 
        g_{j_0, k_0} ,
	& \quad \begin{array}{c} -n_1 \le j_0 \le n_1, \\ -1 \le j_1 \le 1. \end{array}
  \end{alignedt}
\end{equation}
By design,
computing $\widehat{f}_{j_0,j_1}$ for $-n_1 \le j_0 \le n_1$ and $-1 \le j_1 \le
1$ is equivalent to computing $\widehat{f}_j$ for $-n \le j \le n$.

\subsection{Complexity}

If we compute $\widehat{f}_{j_0, j_1}$ in two steps according
to the equations given above, then the number of multiply/adds is 
\[
    N_1^2 N_0 + N N_0 = N ( N_1 + N_0 ) .
\]
On the other hand, the one-step algorithm given by \eqref{1} requires $N^2$ multiply/adds.
Hence, the two-step algorithm beats the one-step algorithm by a factor of
\[
    \frac{N^2}{N ( N_1 + N_0 )} = \frac{N}{N_1 + N_0}
    \approx N/N_1 = N_0 = 3.
\]

\subsection{Recursive Application}

One can do better by iterating the above two-step algorithm.  From the formula
for $g_{j_0,k_0}$ given in \eqref{5}, we see that $g$ is a discrete Fourier transform of 
a subset of the elements of the vector $\{f_k \; : \; k = -n,\ldots,n \}$ 
obtained by sampling $f$ at a
cadence of one every $N_0$ elements.  And, the coefficient $N_0 \dx \dxi$ in the
exponential equals $N_0 / N = 1/N_1$, which again matches the number of terms in
the sum.   Hence, we can apply the two-step algorithm again to this Fourier
transform.  The second key component of the fast Fourier transform is the
observation that this process can be repeated until the Fourier transform only
involves a sum consisting of a single term.

Let $I_N$ denote the number of multiply/adds needed using the recursive
algorithm to solve a problem of size $N = 3^m$.  Keeping in mind that $N_0 = 3$, we get
\begin{eqnarray*}
    I_N = I_{3^m} & = & 3 I_{3^{m-1}} + 3 \cdot 3^m \\
                  & = & 3 (3 I_{3^{m-2}} + 3 \cdot 3^{m-1}) + 3^{m+1} \\
                  & = & 3^2 I_{3^{m-2}} + 2 \cdot 3^{m+1} \\
		  & \vdots & \\
                  & = & 3^k I_{3^{m-k}} + k \cdot 3^{m+1} \\
		  & \vdots & \\
                  & = & 3^m I_{3^{0}} + m \cdot 3^{m+1} \\
                  & = & 3^m (1 + 3 m) \\
                  & = & N (1 + 3 \log_3 N) .
\end{eqnarray*}
Hence, the recursive variant of the algorithm takes on the order of $N \log_3 N$
operations.

\section{A General Factor-Based Algorithm}

The advantage of fast Fourier transforms, such as the one presented in the
previous section, is that they have order $N \log N$ complexity.  But, they have
disadvantages too.  One disadvantage is the need to apply the basic
two-step computation recursively.  Recursion is fine for computing a Fourier
transform, but our aim is to encode a Fourier transform within an optimization model.  
In such a context, it is far better to use a non-recursive algorithm.

A simple modification to the two-step process described in the previous section
produces a variant of the two-step algorithm that makes a more substantial improvement in
the initial two-step computation than what we obtained before.  The idea is to
factor $N$ into a pair of factors with each factor close to the square-root of $N$ rather than
into $3$ and $N/3$.  Indeed, in this section, we assume, as before, that $N$ can be 
factored into
\[
    N = N_0 N_1
\]
but we do not assume that $N_0 = 3$.  In fact, we prefer to have $N_0 \approx
N_1$.  As before, we assume that $N = 2 n + 1$ is odd and therefore that both $N_0$ and
$N_1$ are odd:
\[
    N_0 = 2 n_0 + 1 \quad \text{ and } \quad N_1 = 2 n_1 + 1 .
\]

At the same time, we will now assume that the number of points in the
discretization of the Fourier transform does not necessarily match the number of
points in the discretization of the function itself.  In many real-world
examples, the ``resolution'' of the one discretization does not need to match
the other and artificially enforcing such a match invariably results in a slower
algorithm.  So, suppose that the discrete Fourier transform has the form
\begin{equation} \label{6}
    \widehat{f}_j = \sum_{k=-n}^n e^{2 \pi i k \dx j \dxi} f_k \dx, 
    		\qquad -m \le j \le m,
\end{equation}
and let $M = 2 m + 1$ denote the number of elements in the discretized
transform.  Again, $M$ is odd and therefore we factor it into a product
$M = M_0 M_1$ of two odd factors:
\[
    M_0 = 2 m_0 + 1 \quad \text{ and } \quad M_1 = 2 m_1 + 1 .
\]

If we now decompose our sequencing indices $k$ and $j$ into
\[
    k = N_0 k_1 + k_0 \quad \text{ and } \quad j = M_0 j_1 + j_0,
\]
we get
\begin{eqnarray*} 
&&
    \widehat{f}{j_0,j_1} \\
	    &&
	    \qquad
	= 
	\sum_{k_0=-n_0}^{n_0}
	\sum_{k_1=-n_1}^{n_1} 
	e^{2 \pi i N_0 k_1 \dx M_0 j_1 \dxi} \;
	e^{2 \pi i N_0 k_1 \dx j_0 \dxi} \;
	e^{2 \pi i k_0 \dx (M_0 j_1 + j_0) \dxi} \; \\
		&& \hspace*{3in} \cdot
	f_{k_0,k_1} \dx .
\end{eqnarray*}
As before, we need to assume that the first exponential factor evaluates to
one.  To make that happen, we assume that $N_0 M_0 \dx \dxi$ is an integer.  In
real-world problems, there is generally substantial freedom in the choice of
each of these four factors and therefore guaranteeing that the product is an
integer is generally not a restriction.
With that first exponential factor out of the way, we can again write down a
two-step algorithm
\[
  \begin{alignedt}
    g_{j_0, k_0} 
	& \; = \; &
	\ds \sum_{k_1=-n_1}^{n_1} 
	e^{2 \pi i N_0 k_1 \dx j_0 \dxi}  \;
	f_{k_0,k_1} \dx , 
	& \qquad \begin{array}{c} -m_0 \le j_0 \le m_0, \\
		              -n_0 \le k_0 \le n_0,
			       \end{array}
	\\[0.2in]
    \widehat{f}_{j_0, j_1} 
	& \; = \; &
	\ds \sum_{k_0=-n_0}^{n_0} 
	e^{2 \pi i k_0 \dx (M_0 j_1 + j_0) \dxi} 
        g_{j_0, k_0} ,
	& \qquad \begin{array}{c} -m_0 \le j_0 \le m_0 \\ 
			      -m_1 \le j_1 \le m_1.
			      \end{array}
  \end{alignedt}
\]

\subsection{Complexity}

The number of multiply/adds required for this two-step algorithm is
\[
    NM_0 + MN_0 = MN\left( \frac{1}{M_1} + \frac{1}{N_1} \right).
\]
If $M \approx N$ and $M_1 \approx N_1 \approx \sqrt{N}$, the complexity
simplifies to
\[
    2 N \sqrt{N} .
\]
Compared to the one-step algorithm, which takes $N^2$ multiply/adds, this
two-step algorithm gives an improvement of a factor of $\sqrt{N}/2$.   This
first-iteration improvement is much better than the factor of $3$ improvement
from the first iteration of the recursive algorithm of the previous section.
Also, if $M$ is much smaller than $N$, we get further
improvement over the full $N \times N$ case.

Of course, if $M_0, M_1, N_0$, and $N_1$ can be further factored, then this
two-step algorithm can be extended in the same manner as was employed for the
algorithm of the previous section successively factoring $M$ and $N$ until it is
reduced to prime factors.  But, our eventual aim in this paper is to
embed these algorithms into an optimization algorithm and so we will focus our
attention in this paper just on two-step algorithms and not their recursive
application.

\section{Fourier Transforms in 2D}

Many real-world optimization problems, and in particular the one to be
discussed in Section \ref{sec_appl}, involve Fourier transforms in more than
one dimension.  It turns out that the core idea in the algorithms presented
above, replacing a one-step computation with a two-step equivalent, presents
itself in this higher-dimensional context as well \cite{SPSV07}.

Consider a two-dimensional Fourier transform
\[
    \widehat{f}(\xi, \eta) 
	=
	\iint e^{2 \pi i (x \xi + y \eta)} f(x,y) dy dx 
\]
and its discrete approximation
\[
    \widehat{f}_{j_1, j_2}
	=
	\sum_{k_1=-n}^{n} \sum_{k_2=-n}^{n} 
        e^{2 \pi i (x_{k_1} \xi_{j_1} + y_{k_2} \eta_{j_2})} f_{k_1, k_2} \dy \dx , 
		\qquad -m \le j_1, j_2 \le m,
\]
where 
\[
  \begin{alignedt}
    x_k & = k \dx, & \qquad -n \le k \le n, \\
    y_k & = k \dy, & \qquad -n \le k \le n, \\
    \xi_j & = j \dxi, & \qquad -m \le j \le m, \\
    \eta_j & = j \deta, & \qquad -m \le j \le m, \\
    f_{k_1, k_2} & = f(x_{k_1}, y_{k_2}), & \qquad -n \le k_1, k_2 \le n \\
    \widehat{f}_{j_1, j_2} & = \widehat{f}(\xi_{j_1}, \eta_{j_2}), 
    			& \qquad -m \le j_1, j_2 \le m.
  \end{alignedt}
\]

Performing the calculation in the obvious way requires $M^2 N^2$ complex additions and
a similar number of multiplies.  However, we can factor the exponential into the
product of two exponentials and break the process into two steps:
\[
  \begin{alignedt}
    g_{j_1, k_2} 
        & = \ds \sum_{k_1=-n}^{n} e^{2 \pi i x_{k_1} \xi_{j_1}} f_{k_1, k_2} \dx , 
	& \qquad -m \le j_1 \le m, -n \le k_2 \le n, \\[0.2in]
    \widehat{f}_{j_1, j_2} 
	& = \ds \sum_{k_2=-n}^{n} e^{2 \pi i y_{k_2} \eta_{j_2}} g_{j_1, k_2} \dy , 
	& \qquad -m \le j_1, j_2 \le m,
  \end{alignedt}
\]
It is clear that, in this context, the two-step approach is simply to break up
the two-dimensional integral into a nested pair of one-dimensional integrals.
Formulated this way, the calculation requires only $MN^2 + M^2N$ complex
additions (and a similar number of multiplications).

The real-world example we shall discuss shortly involves a two-dimensional
Fourier transform.  Given that the idea behind speeding up a one-dimensional
Fourier transform is to reformulate it as a two-dimensional transform and then
applying the two-step speed up trick of the two-dimensional transform, we shall for
the rest of the paper restrict our attention to problems that are two
dimensional.

\section{Exploiting Symmetry}

Before discussing real-world examples and associated computational results, it
is helpful to make one more simplifying assumption.  If we
assume that $f$ is invariant under reflection about both the $x$ and $y$
axes, i.e., $f(-x,y) = f(x,y)$ and $f(x,-y) = f(x,y)$ for all $x$ and $y$, then
the transform has this same symmetry and is in fact real-valued.  In this case,
it is simpler to use an even number of grid-points ($N = 2n$ and $M = 2m$) 
rather than an odd number and
write the straightforward algorithm for the two-dimensional discrete Fourier transform as
\begin{equation} \label{13}
  \begin{alignedt}
    \widehat{f}_{j_1, j_2} 
	& = \ds 4 
	\sum_{k_1=1}^{n} \sum_{k_2=1}^{n} 
	\cos(2 \pi x_{k_1} \xi_{j_1}) \cos(2 \pi y_{k_2} \eta_{j_2}) 
	f_{k_1, k_2} \dy \dx, 
	& \qquad 1 \le j_1, j_2 \le m,
  \end{alignedt}
\end{equation}
where
\[
  \begin{alignedt}
    x_k & = & (k-1/2) \dx, & \qquad 1 \le k \le n, \\[0.1in]
    y_k & = & (k-1/2) \dy, & \qquad 1 \le k \le n, \\[0.1in]
    \xi_j & = & (j-1/2) \dxi, & \qquad 1 \le j \le m, \\[0.1in]
    \eta_j & = & (j-1/2) \deta, & \qquad 1 \le j \le m, \\[0.1in]
    f_{k_1, k_2} & = & f(x_{k_1}, y_{k_2}), & \qquad 1 \le k_1, k_2 \le n \\[0.1in]
    \widehat{f}_{j_1, j_2} & \approx & \widehat{f}(\xi_{j_1}, \eta_{j_2}), 
    			& \qquad 1 \le j_1, j_2 \le m. \\~
  \end{alignedt}
\]
The two-step algorithm then takes the following form:
\[
  \begin{alignedt}
    g_{j_1, k_2} 
        & = \ds 2 \sum_{k_1=1}^{n} \cos(2 \pi x_{k_1} \xi_{j_1}) f_{k_1, k_2} \dx , 
	& \qquad 1 \le j_1 \le m, 1 \le k_2 \le n, \\[0.2in]
    \widehat{f}_{j_1, j_2} 
	& = \ds 2 \sum_{k_2=1}^{n} \cos(2 \pi y_{k_2} \eta_{j_2}) g_{j_1, k_2} \dy , 
	& \qquad 1 \le j_1, j_2 \le m,
  \end{alignedt}
\]

\subsection{Complexity}

The complexity of the straightforward one-step algorithm is $m^2 n^2$ 
and the complexity of the two-step algorithm is $mn^2 + m^2n$.  Since $m = M/2$
and $n = N/2$, we see that by exploiting symmetry 
the straightforward algorithm gets speeded up by
a factor of $16$ and the two-step algorithm gets speeded up by a factor of
$8$.   But, the improvement is better than that as complex arithmetic has also
been replaced by real arithmetic.  One complex add is the same as two real adds
and one complex multiply is equivalent to four real multiplies and two real
adds.   Hence, complex arithmetic is about four times more computationally
expensive than real arithmetic.

\section{Matrix Notation}
As Fourier transforms are linear operators it is instructive to express our
algorithms in matrix/vector notation.  In this section, we shall do this for the
two-dimensional Fourier transform.  To this end, 
let $F$ denote the $n \times n$ matrix with elements $f_{k_1,k_2}$, 
let $G$ denote the $m \times n$ matrix with elements $g_{j_1,k_2}$,
let $\widehat{F}$ denote the $m \times m$ matrix with elements $\widehat{f}_{j_1,j_2}$,
and let $K$ denote the $m \times n$ Fourier kernel matrix whose elements are
\[
    \kappa_{j_1,k_2} = \cos(2 \pi x_{k_1} \xi_{j_1}) \dx . 
\]
For notational simplicity,
assume that the discretization in $y$ is the same as it is in $x$, 
	i.e., $\dx = \dy$, 
and that the discretization in $\eta$ is the same as it is in $\xi$,
	i.e., $\deta = \dxi$.
Then, the two-dimensional Fourier transform $\widehat{F}$ can be written simply
as
\[
    \widehat{F} = K F K^T 
\]
and the computation of the transform in two steps is just the statement that the
two matrix multiplications can, and should, be done separately:
\[
  \begin{alignedt}
     G & = K F \\[0.1in]
     \widehat{F} & = G K^T .
  \end{alignedt}
\]

When linear expressions are passed to a linear programming code, the variables
are passed as a vector and the constraints are expressed in terms of a matrix of
coefficients times this vector.   The matrix $F$ above represents the variables in
the optimization problem.  If we let $f_k$, $k=1,\ldots,n$ denote the $n$
columns of this matrix, i.e., $F = [ f_{1} \; f_{2} \; \cdots \; f_n ]$,
then we can list the elements in column-by-column order to make a column vector 
(of length $n^2$):
\[
    \text{vec}(F) =
    \left[
    \begin{array}{c}
	f_{1} \\
	f_{2} \\
	\vdots \\
	f_n
    \end{array}
    \right] .
\]
Similarly, we can list the elements of $G$ and $\widehat{F}$ in column vectors too:
\[
    \text{vec}(G) =
    \left[
    \begin{array}{c}
	g_{1} \\
	g_{2} \\
	\vdots \\
	g_n
    \end{array}
    \right] 
    \qquad
    \text{and}
    \qquad
    \text{vec}(\widehat{F}) =
    \left[
    \begin{array}{c}
	\widehat{f}_{1} \\
	\widehat{f}_{2} \\
	\vdots \\
	\widehat{f}_m
    \end{array}
    \right] 
    .
\]
It is straightforward to check that
\[
    \text{vec}(G) 
	=
	\left[
	\begin{array}{cccc}
	    K \\
	    & K \\
	    && \ddots \\
	    &&& K 
	\end{array}
	\right]
        \text{vec}(F) 
\]
and that
\[
    \text{vec}(\widehat{F}) 
	=
	\left[
	\begin{array}{cccc}
	    \kappa_{1,1}I \; & \kappa_{1,2}I \; & \cdots \; & \kappa_{1,n}I \\
	    \kappa_{2,1}I \; & \kappa_{2,2}I \; & \cdots \; & \kappa_{2,n}I \\
	    \vdots        \; & \vdots        \; & \ddots \; & \vdots \\
	    \kappa_{m,1}I \; & \kappa_{m,2}I \; & \cdots \; & \kappa_{m,n}I \\
	\end{array}
	\right]
        \text{vec}(G) ,
\]
where $I$ denotes an $m \times m$ identity matrix.

The matrices in these two formulae are sparse:  the first is block diagonal and
the second is built from identity matrices.   Passing the constraints to a
solver as these two sets of constraints introduces new variables and more
constraints, but the constraints are very sparse.   Alternatively, if we were to
express $\text{vec}({\widehat{F}})$ directly in terms of $\text{vec}(F)$, these
two sparse matrices would be multiplied together and a dense coefficient matrix
would be passed to the solver.  It is often the case that optimization problems
expressed in terms of sparse matrices solve much faster than equivalent
formulations involving dense matrices even when the latter involves fewer
variables and/or constraints (see, e.g., \cite{Van90e}).

\section{A Real-World Example: High-Contrast Imaging} \label{sec_appl}

Given the large number of planets discovered over the past decade by so-called
``indirect'' detection methods, there is great interest in building a special
purpose telescope capable of imaging a very faint planet very close to
its much brighter host star.  This is a problem in {\em high-contrast
imaging}.   It is made difficult by the fact that light is a wave and therefore
point sources, like the star and the much fainter planet, produce not just 
single points of light in the image but rather diffraction patterns---most of
the light lands where ray-optics suggests it will but some of the light lands
nearby but not exactly at this point.   In a
conventional telescope, the ``wings'' of the diffraction pattern produced by the
star are many
orders of magnitude brighter than any planet would be at the place where the
planet might be.   Hence, the starlight outshines the planet and makes the
planet impossible to detect.   But, it is possible to customize the diffraction
pattern by designing an appropriate filter, or a mask, to put on the front of the telescope.
While it is impossible to concentrate all of the starlight at the central
point---to do so would violate the uncertainty principle---it is
possible to control it in such a way that there is a very dark patch very close
to the central spot.

Suppose that we place a filter over the opening of a telescope with the
property that the transmissivity of the filter varies from place to place over the
surface of the filter.   Let $f(x,y)$ denote the transmissivity at location
$(x,y)$ on the surface of the filter ($(0,0)$ denotes the center of the filter).
It turns out that the {\em electromagnetic field} in the image plane of such a telescope
associated with a single point on-axis source (the star) is proportional to the Fourier
transform of the filter function $f$.  
Choosing units
in such a way that the telescope's opening has a diameter of one,
the Fourier transform can be written as
\begin{equation} \label{22}
    \widehat{f}(\xi, \eta) = \int_{-1/2}^{1/2} \int_{-1/2}^{1/2}
    e^{2 \pi i (x \xi + y \eta)} f(x,y) dy dx .
\end{equation}
The {\em intensity} of the light in the image
is proportional to the magnitude squared of the electromagnetic field.  

Assuming that the underlying telescope has a circular opening of diameter one,
we impose the following constraint on the function $f$:
\[
    f(x,y) = 0 \quad \text{ for } \quad x^2 + y^2 > (1/2)^2 .
\]

As often happens in real-world problems, there are multiple competing goals.
We wish to maximize the amount of light that passes through the filter and
at the same time minimize the amount of light that lands within a dark zone
$\mathcal{D}$ of the image plane.  If too much light lands in the dark zone, the
telescope will fail to detect the planets it is designed to find.  Hence, this
latter objective is usually formulated as a constraint.  This leads to the
following optimization problem:
\begin{equation} \label{8}
  \begin{array}{ll}
    \text{maximize }   & \ds \iint f(x,y) dy dx \\
    \text{subject to } &
      \begin{alignedt}
          & \left|\widehat{f}(\xi, \eta)\right|^2 & \le \varepsilon, 
		  & \quad (\xi, \eta) \in \mathcal{D}, \\
	  & f(x,y)                 & = 0,          
		  & \quad x^2 + y^2 > (1/2)^2,         \\
	  0 \le & f(x,y) & \le 1, & \quad \text{ for all } x,y  .
      \end{alignedt}
  \end{array}
\end{equation}
Here, $\varepsilon$ is a small positive constant representing the maximum level of
brightness of the starlight in the dark zone.  Without imposing further symmetry
constraints on the function $f$, the Fourier transform $\widehat{f}$ is complex
valued.  Hence this optimization problem has a
linear objective function and both linear constraints and convex quadratic inequality constraints.   
Hence, a discretized version can be solved (to a global optimum) using, say, interior-point methods.

Assuming that the filter can be symmetric with respect to reflection about both
axes (in real-world examples, this is often---but not always---possible; see
		\cite{CVK11b} for several examples),
the Fourier transform can be written as
\[
    \widehat{f}(\xi, \eta) = 4 \int_0^{1/2} \int_0^{1/2}
    \cos(2 \pi x \xi ) \cos(2 \pi y \eta) f(x,y) dy dx .
\]
In this case, the Fourier transform is real and so the convex quadratic inequality
constraint in \eqref{8} can be replaced with a pair of inequalities,
\[
    -\sqrt{\varepsilon} \le \widehat{f}(\xi, \eta) \le \sqrt{\varepsilon} ,
\]
making the problem an infinite dimensional linear programming problem.

Figure \ref{fig1} shows an \ampl\ model formulation of this problem expressed in
the straightforward one-step manner.   Figure \ref{fig2}, on the other hand,
shows an \ampl\ model for the same problem but with the Fourier transform
expressed as a pair of transforms---the so-called two-step process.

\begin{figure*}
\begin{center}
\scriptsize
\begin{verbatim}
param pi := 4*atan(1);
param rho0 := 4;
param rho1 := 20;

param n := 150;          # discretization parameter
param dx := 1/(2*n);
param dy := dx;
set Xs := setof {j in 0.5..n-0.5 by 1} j/(2*n);
set Ys := setof {j in 0.5..n-0.5 by 1} j/(2*n);
set Pupil := setof {x in Xs, y in Ys: x^2+y^2 < 0.25} (x,y);

var f {x in Xs, y in Ys: x^2 + y^2 < 0.25} >= 0, <= 1, := 0.5;

param m := 35;          # discretization parameter
set Xis := setof {j in 0..m} j*rho1/m;
set Etas := setof {j in 0..m} j*rho1/m;
set DarkHole := setof {xi in Xis, eta in Etas: 
                        xi^2+eta^2>=rho0^2 && 
                        xi^2+eta^2<=rho1^2 && 
                        eta <= xi } (xi,eta);

var fhat {xi in Xis, eta in Etas} = 
    4*sum {(x,y) in Pupil} f[x,y]*cos(2*pi*x*xi)*cos(2*pi*y*eta)*dx*dy;

var area = sum {(x,y) in Pupil} f[x,y]*dx*dy;

maximize throughput: area;

subject to sidelobe_pos {(xi,eta) in DarkHole}: fhat[xi,eta] <= 10^(-5)*fhat[0,0];
subject to sidelobe_neg {(xi,eta) in DarkHole}: -10^(-5)*fhat[0,0] <= fhat[xi,eta];

solve;
\end{verbatim}
\end{center}
\caption{\ampl\ model for discretized version of problem \eqref{8} assuming that
the mask is symmetric about the $x$ and $y$ axes.  The dark zone $\mathcal D$ is
a pair of sectors of an annulus with inner radius $4$ and outer radius $20$.
The optimal solution is shown in Figure \ref{fig3}.}
\label{fig1}
\end{figure*}

\begin{figure*}
\begin{center}
\scriptsize
\begin{verbatim}
param pi := 4*atan(1);
param rho0 := 4;
param rho1 := 20;

param n := 1000;                # discretization parameter
param dx := 1/(2*n);
param dy := dx;
set Xs := setof {j in 0.5..n-0.5 by 1} j/(2*n);
set Ys := setof {j in 0.5..n-0.5 by 1} j/(2*n);
set Pupil := setof {x in Xs, y in Ys: x^2+y^2 < 0.25} (x,y);

var f {x in Xs, y in Ys: x^2 + y^2 < 0.25} >= 0, <= 1, := 0.5;

param m := 35;                # discretization parameter
set Xis := setof {j in 0..m} j*rho1/m;
set Etas := setof {j in 0..m} j*rho1/m;
set DarkHole := setof {xi in Xis, eta in Etas: 
                        xi^2+eta^2>=rho0^2 && 
                        xi^2+eta^2<=rho1^2 && 
                        eta <= xi } (xi,eta);

var g {xi in Xis, y in Ys};
var fhat {xi in Xis, eta in Etas};

var area = sum {(x,y) in Pupil} f[x,y]*dx*dy;

maximize throughput: area;

subject to g_def {xi in Xis, y in Ys}:
    g[xi,y] = 2*sum {x in Xs: (x,y) in Pupil} 
        f[x,y]*cos(2*pi*x*xi)*dx;

subject to fhat_def {xi in Xis, eta in Etas}:
    fhat[xi,eta] = 2*sum {y in Ys} 
        g[xi,y]*cos(2*pi*y*eta)*dy;

subject to sidelobe_pos {(xi,eta) in DarkHole}:  fhat[xi,eta] <= 10^(-5)*fhat[0,0];
subject to sidelobe_neg {(xi,eta) in DarkHole}: -10^(-5)*fhat[0,0] <= fhat[xi,eta];

solve;
\end{verbatim}
\end{center}
\caption{\ampl\ model reformulated to exploit the two-step algorithm.  The
	optimal solution is shown in Figure \ref{fig4}.}
\label{fig2}
\end{figure*}

\begin{figure*}
\begin{center}
\includegraphics[width=2.3in]{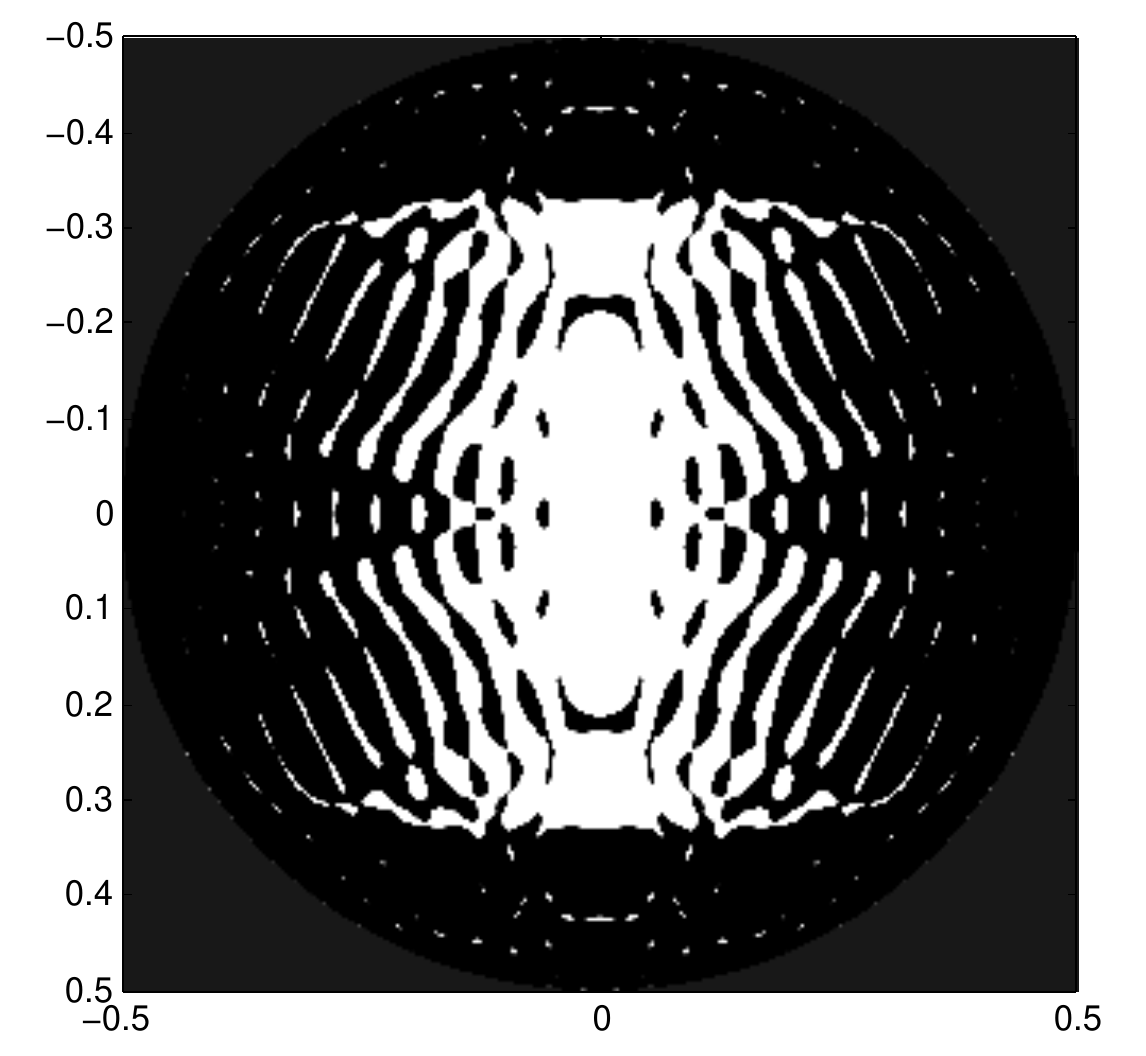}
\includegraphics[width=2.3in]{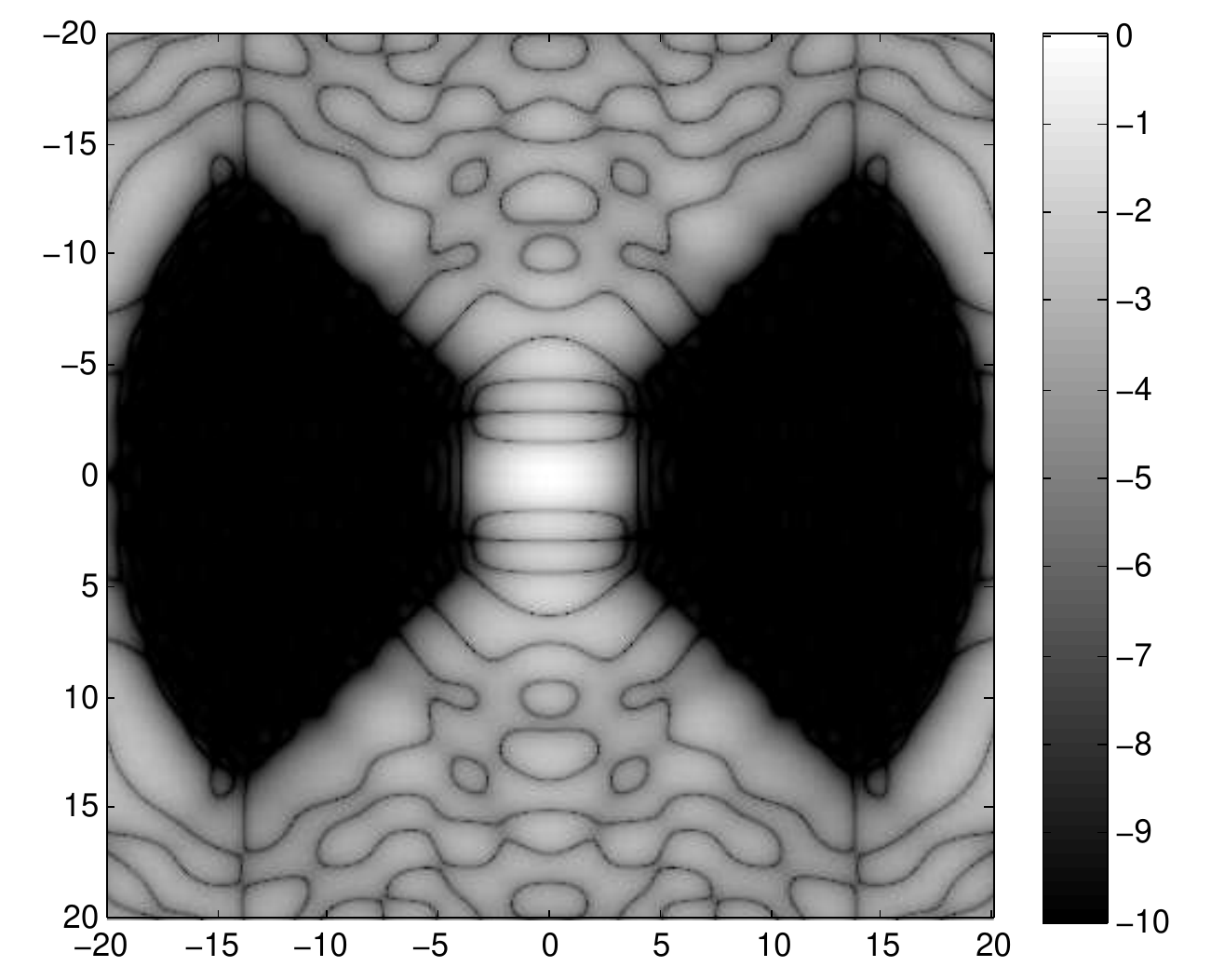}
\end{center}
\caption{The optimal filter from the one-step model shown in Figure \ref{fig1}, 
	which turns out to be purely opaque and transparent
	(i.e., a mask), and a logarithmic plot of the star's image.}
	\label{fig3}
\end{figure*}

\begin{figure*}
\begin{center}
\includegraphics[width=2.3in]{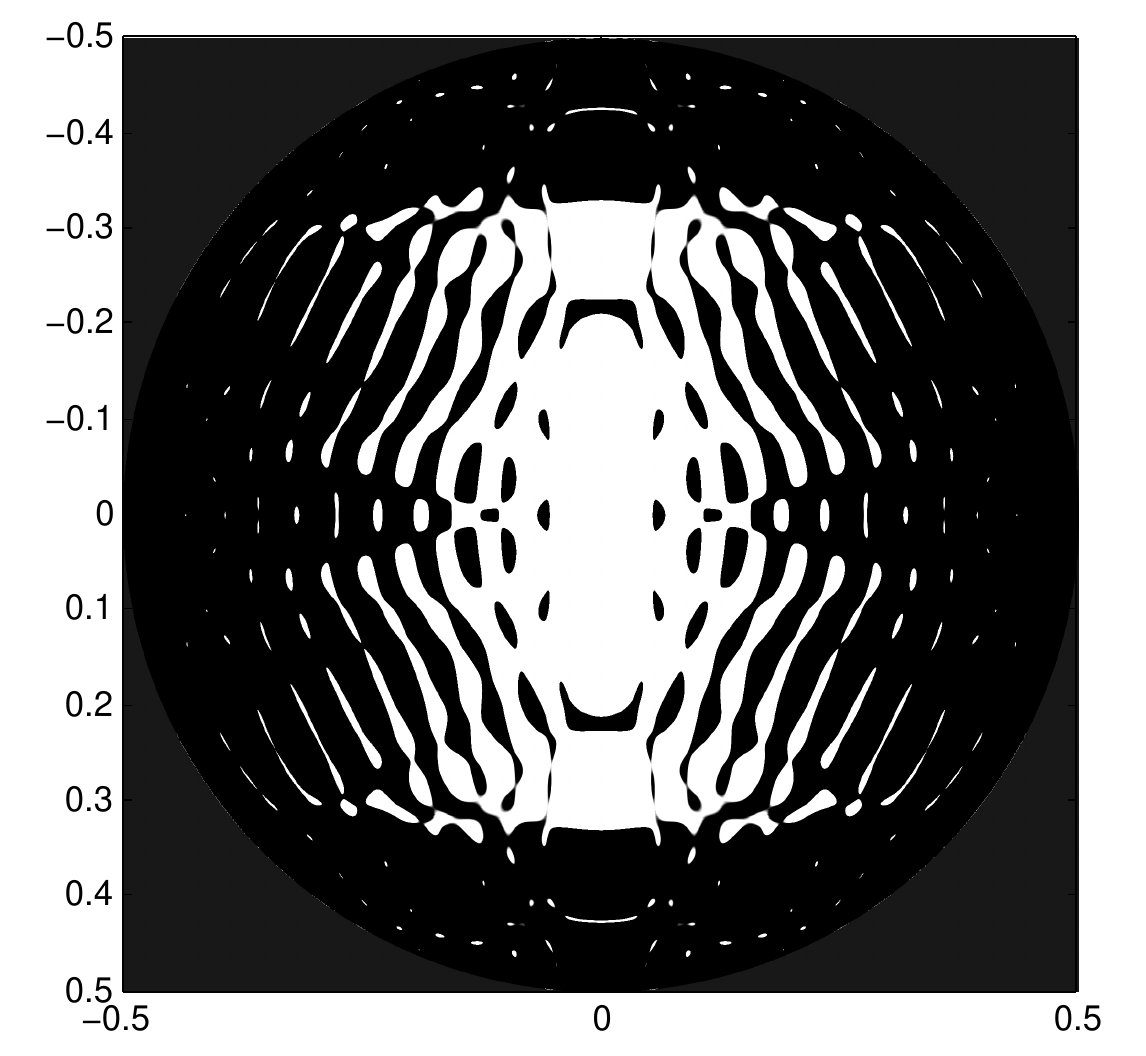}
\includegraphics[width=2.3in]{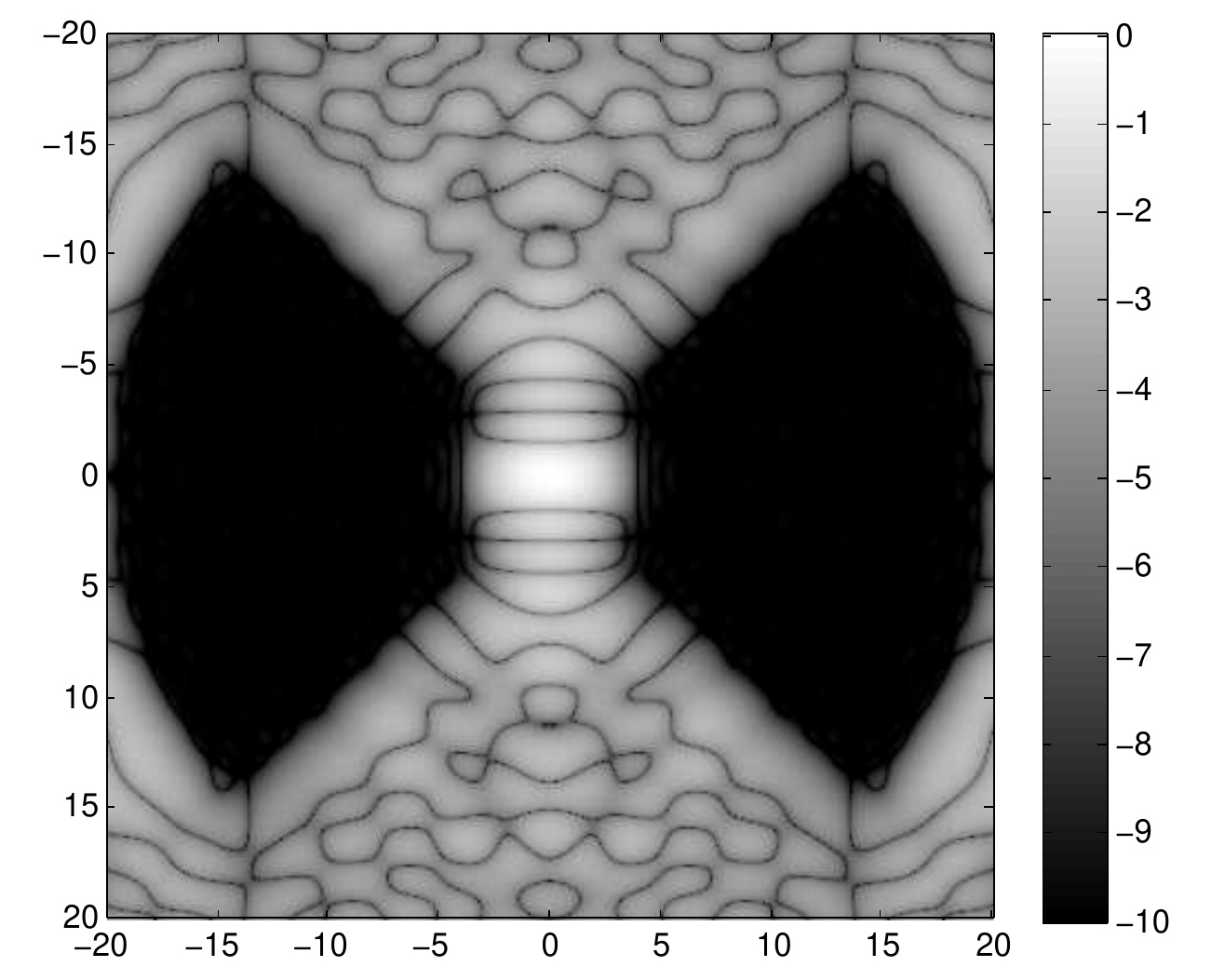}
\end{center}
\caption{The optimal filter from the two-step model shown in Figure \ref{fig2} 
	and a logarithmic plot of the star's image.}
	\label{fig4}
\end{figure*}

\begin{figure*}
\begin{center}
\includegraphics[width=2in]{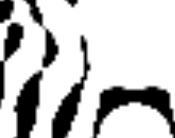}
\hfill
\includegraphics[width=2in]{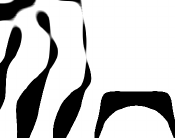}
\end{center}
\caption{Close up of the two masks to compare resolution.}
	\label{fig5}
\end{figure*}

As Figures \ref{fig3}, \ref{fig4}, and \ref{fig5} show,
the optimal solution for the two models are, of course, essentially the same
except for the improved resolution in the two-step version provided by 
a larger value for $n$ ($n=1000$ vs. $n=150$).  Using \loqo\ \cite{Van94f} as the
interior-point method to solve the problems, both versions solve in
a few hours on a modern computer.  It is possible to solve even larger
instances, say $n=2000$, if one is willing to wait a day or so for a solution.
Ultimately, higher resolution is actually important
because manufacturing these masks involves replacing the pixellated mask with a
spline-fitted smooth approximation and it is important to get this approximation
correct.

Table \ref{tab1} summarizes problem statistics for the two versions of the
model as well as a few other size choices.
Table \ref{tab2} summarizes solution statistics for these same problems.  These
problems were run as a single thread on a GNU/Linux 
(Red Hat Enterprise Linux Server release 5.7) x86\_64 server with dual
Xeon X5460s cpus (3.16 GHz with 4 cores each), 32 GB of RAM and a 6.1 MB cache.

\begin{table}
\caption{Comparison between a few sizes of the one-step model shown in Figure \ref{fig1}
	and a few sizes of the two-step model shown in Figure \ref{fig2}.
	The column labeled {\tt nonzeros} reports the number of nonzeros in the
	constraint matrix of the linear programming problem and the column 
	{\tt arith. ops.}
	The One-Step-250x35
	problem is too large to solve by \loqo, which is compiled for a 32-bit architecture
	operating system.  
	}
	\label{tab1}
\begin{center}
\begin{tabular}{lrr|rrrrrrrr}
\hline
  Model    & $n$ & $m$ & constraints & variables &   nonzeros &    arith. ops. \\ \hline
  One step & 150 & 35 &    976 &  17,672 & 17,247,872 & 17,196,541,336 \\
  One step & 250 & 35 &      * &       * &          * &              * \\ \hline
  Two step & 150 & 35 &  7,672 &  24,368 &    839,240 &  3,972,909,664 \\
  Two step & 500 & 35 & 20,272 & 215,660 &  7,738,352 & 11,854,305,444 \\
  Two step &1000 & 35 & 38,272 & 822,715 & 29,610,332 & 23,532,807,719 \\
\hline
\end{tabular}
\end{center}
\vspace*{0.3in}
\end{table}

\begin{table}
\caption{Hardware-specific performance comparison data.  
	The results shown here were obtained using the
	default value for all of \loqo's tunable parameters.  It is possible to
	reduce the iteration counts to about 100 or less on all the problems by
	increasing the value of the {\tt epsdiag} parameter to about {\tt 1e-9}.}
	\label{tab2}
\begin{center}
\begin{tabular}{lrr|rrrrrrrr}
\hline
  Model    & $n$ & $m$ & iterations & primal objective & dual objective &  cpu time (sec) \\
\hline
One step &  150 & 35 &  54 & 0.05374227247 & 0.05374228041 &  1380 \\
One step &  250 & 35 &   * &             * &             * &     * \\ \hline
Two step &  150 & 35 & 185 & 0.05374233071 & 0.05374236091 &  1064 \\
Two step &  500 & 35 & 187 & 0.05395622255 & 0.05395623990 &  4922 \\
Two step & 1000 & 35 & 444 & 0.05394366337 & 0.05394369256 & 26060 \\
\hline
\end{tabular}
\end{center}
\vspace*{0.3in}
\end{table}


Real telescopes have opennings that are generally not just open
unobstructed disks but, rather, typically have central obstructions supported by
spiders.  It is easy to extend the ideas presented here to handle such
situations;  see \cite{CVK11b}.

As explained in earlier sections, the two-step algorithm applied to a
one-dimensional Fourier transform effectively makes a two-dimensional
representation of the problem and applies the same two-step algorithm that we
have used for two-dimensional Fourier transforms.  It is natural, therefore to
consider whether we can get more efficiency gains by applying the two-step
algorithm to each of the iterated one-dimensional Fourier transforms that make
up the two-step algorithm for the two-dimensional Fourier transform.
We leave such investigations for future work.


\begin{figure*}
\begin{center}
\includegraphics[width=2.3in]{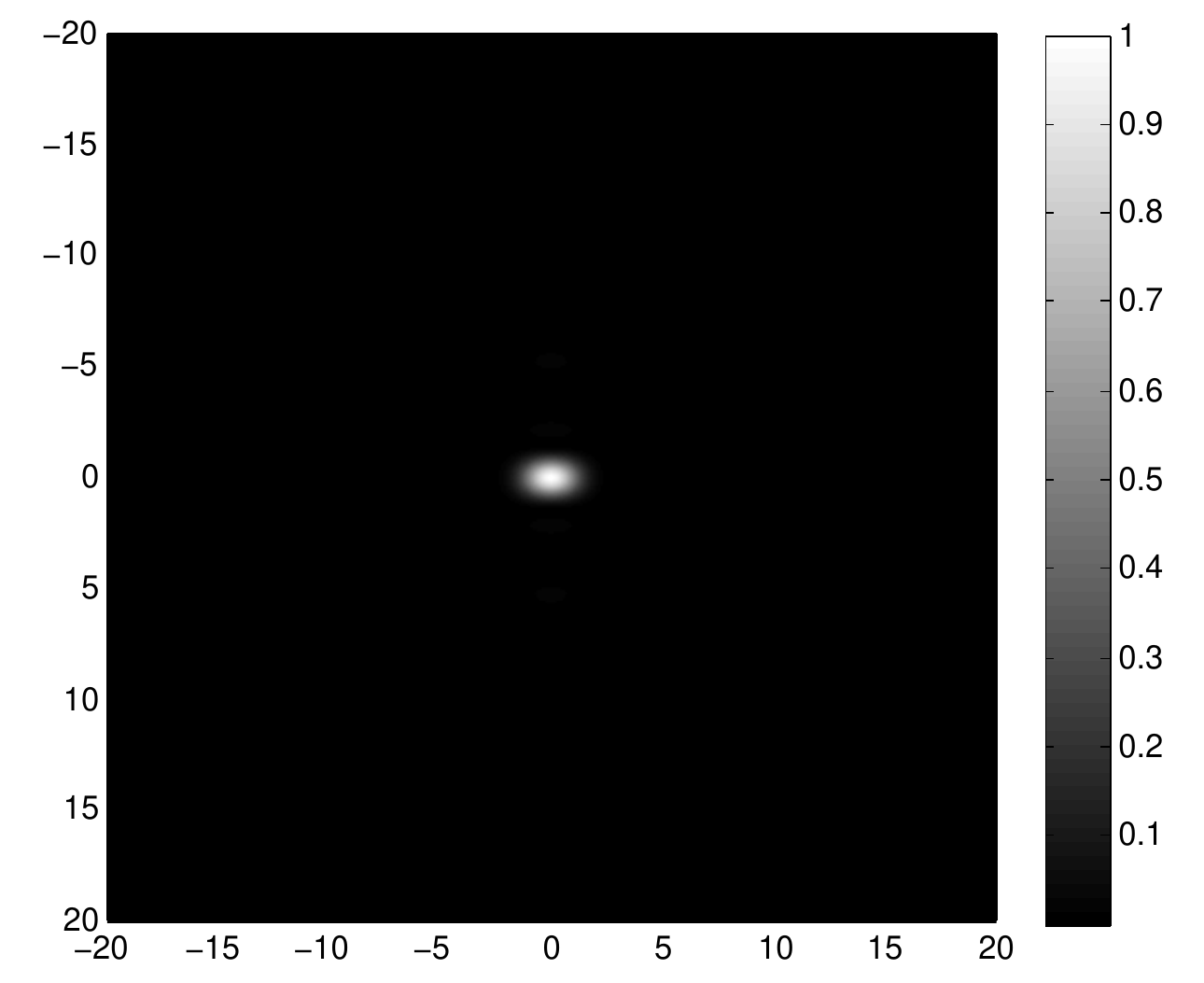}
\includegraphics[width=2.3in]{TwoSectorDarkHole_Sparse1000_psf.pdf}
\end{center}
\caption{Logarithmic stretches are useful but can be misleading. 
	{\em Left:} the image of the star shown in a linear stretch.
	{\em Right:} the same image shown in a logarithmic stretch.}
	\label{fig6}
\end{figure*}

%
%

\begin{acknowledgements}
I would like to thank N. Jeremy Kasdin, Alexis Carlotti, and all the members of
the Princeton High-Contrast Imaging Lab for many stimulating discussions.
\end{acknowledgements}

\bibliographystyle{spmpsci}      
\bibliography{../lib/refs}   


\end{document}